\documentclass[12pt]{article}
\newlength{\wdth}

\usepackage{hyperref}

\newif\ifhide
\hidefalse

\def\P{\mathsf{P}}
\def\E{\mathsf{E}}
\def\F{{\cal F}}

\newcommand{\R}        {{{\rm I \hskip -2pt R}}}

\newtheorem{lemma}{Lemma}

\newtheorem{theorem}{Theorem}
\newtheorem{Proposition}{Proposition}

\newtheorem{Corollary}{Corollary}

\newtheorem{Remark}{Remark}

\usepackage{amsmath}
\usepackage{amssymb}
\usepackage{color}
\usepackage{hyperref}

\newif\ifproofs
\proofsfalse

\newif\ifp
\pfalse

\begin{document}
\title{{\normalsize\tt\hfill\jobname.tex}\\
On higher order moments and rates of convergence for SDEs with  
switching 
}
\author{Alexander Veretennikov\footnote{
Institute for Information Transmission Problems, Moscow, Russia;  email: {ayv@iitp.ru}}\;\footnote{This study was funded by the Foundation for the Advancement of Theoretical Physics and Mathematics ``BASIS''.}
}

\maketitle       

%\tableofcontents

\begin{abstract}
Second order recurrence of a $d$-dimensional diffusion with an additive Wiener process,  with switching,  and with one recurrent and one transient regime and constant switching intensities is established under suitable conditions. The approach is based on embedded Markov chains and a priori bounds for the moments of $X_t$ at moments of jumps of the discrete component, as well as on certain simple martingale properties.

%{\color{red}Soslat'sia na Belopolskuyu o sistemah, tam nuzhny diffuzii s perekliucheniem} section 2.2

\medskip


\medskip


\end{abstract}
\section{Introduction}
We are interested in the second order recurrence, or 2-recurrence for the process $(X_t,Z_t)$ with a continuous component $X$ and discrete one $Z$ described by the stochastic differential equation  in $\R^d$
\begin{align}\label{sde}
dX_{t} =b(X_{t}, Z_t)\, dt+ 
dW_{t}, \quad t\ge 0, 
\quad X_{0} =x, \; Z_0=z, 
\end{align}
for the component $X$, while $Z_t$ is a continuous-time  Markov process on the state space $S= \{0,1\}$ with constant positive intensities of respective transitions  $\lambda_{01} =: \lambda_-, \, \& \, \lambda_{10} =: \lambda_+$;  the trajectories of $Z$ are assumed to be c\`adl\`ag; the jumps for $Z$ are independent of $W$ and of the trajectory of the component $X$. 
Positive recurrence for such systems was recently studied in \cite{Ver21, Ver21b}, and the history and other references are briefly commented in the end of this section. The second order recurrence is useful because it provides a certain rate of convergence better than what follows from the mere positive recurrence; of course, it 
will require some additional qualitative assumptions in comparison to the conditions for the positive recurrence. On the other hand, the assumptions required for such recurrence is -- quite naturally -- weaker than what is needed for the exponential convergence (see \cite{Hairer}).
The notations
$$
b(x,0) = b_-(x), \quad b(x,1) = b_+(x), 
$$
will be used. The condition of the boundedness of the function $b$  suffices for the process $(X_t,Z_t)$ to be well-defined as a strong solution; see  \cite{Ver79, Ver21, Ver21b} for some details; see also \cite[section 2.2]{Yana} concerning the links between semi-linear parabolic systems and SDEs with switching. 
The process $(X,Z)$ is a Feller one and strong Markov's.

\medskip

The SDE solution is assumed ergodic under the regime $Z=0$ and transient under $Z=1$. We are looking for sufficient conditions for the second order recurrence of the strong Markov process $(X_t,Z_t)$. 
Such a problem was considered in \cite{Hairer} for the exponentially recurrent case; positive recurrence was studied in \cite{Ver21} and \cite{Ver21b}; for other references see \cite{Anulova}, \cite{Khasminskii12}, \cite{Mao}, \cite{ShaoYuan}, and the references therein. 
Under weak ergodic and transient conditions the setting was earlier investigated in \cite{Ver21} for the case of the constant intensities $\lambda_{0}, \lambda_{1}$ (i.e., not depending on $x$). 
The result from \cite{Ver21} will be used essentially; however, the upper bounds for the second moment of $X$ from \cite{Ver21} do not suffice,  and some kind of similar {\em lower} bounds for the second moments of the diffusion component will be  established under additional assumptions: it turns out that they are also needed for the upper bound of the second order recurrence. 

The paper consists of the sections: Introduction, Main results (section \ref{sec_main}), Auxiliaries I (section \ref{sec_auxI}), Auxiliaries II (section \ref{sec_auxII}), Auxiliaries III (section \ref{sec_auxIII}), and Proof of main result, which is theorem \ref{thm2} (section \ref{sec_mainproof}). 

Sometimes we will drop the indices from the expectations writing $\E$ instead of $\E_{x,z}$ where there is no confusion. Sometimes we will drop just the index $z$ where appropriate.

\section{Main results}\label{sec_main}
In \cite{Ver21} and in \cite{Ver21b} the following result was established (in fact, in \cite{Ver21} for $d=1$ and constant $\lambda_-, \lambda_+$, 
and in \cite{Ver21b} for $d\ge 1$ and for variable coefficients 
$\lambda_-, \lambda_+$). 
\begin{Proposition}[\cite{Ver21}]\label{thm1}
Let 
the drift $b = (b_+,b_-)$ be bounded and Borel measurable, and let there exist $r_-, r_+,M>0$ such that 
\begin{equation}\label{al}
0<{\lambda}_- \wedge{\lambda}_+ 
\le  \lambda_- \vee \lambda_1 < \infty,
\end{equation}
\begin{equation}\label{b}
x b_-(x) \le - r_-, \quad x b_+(x)\le + r_+, \quad \forall \,|x|\ge M, 
\end{equation}
and 
\begin{equation}\label{c1}
2r_- > d \; 
\quad \& \quad \frac{(2r_- -  d)}{\lambda_-} > \frac{(2r_++ d)}{{\lambda}_+}.
\end{equation}
Then the process $(X,Z)$ is positive recurrent; moreover, 
there exists $C>0$ such that for all $M_1$ large enough and all $x \in \mathbb R^d$ and for $z=0,1$
\begin{equation}\label{e3}
\mathsf E_{x,z}\tau_{M_1} 
\le C (x^2 + 1), 
\end{equation} 
where 
$$
\tau_{M_1} := \inf(t\ge 0:\, |X_t|\le M_1). 
$$
Moreover, the process $(X_t,Z_t)$ has a unique invariant measure, and for each nonrandom initial condition $x,z$ there is a convergence to this measure in total variation when $t\to\infty$.
\end{Proposition}

\noindent
For studying the second order recurrence we assume a bit more restrictive condition:
\begin{equation}\label{c2}
\frac{(4r_--(2d+4))}{\lambda_-} > \frac{(4r_++(2d+4))}{\lambda_+}.
\end{equation}
Note that (\ref{c2}) implies (\ref{c1}). 
In fact, for the proof, {\em really for simplicity of the presentation,}  we assume even a bit more, namely,  
\begin{equation}\label{c2a}
\frac{(6r_--(3d+12))}{\lambda_-} > \frac{(6r_++(3d+12))}{\lambda_+}.
\end{equation}
Note that (\ref{c2a}) implies (\ref{c2}). 
Also, it will be assumed the following:
\begin{equation}\label{b2}
x b_-(x) \ge - R_-, \quad x b_+(x)\ge + R_+, \quad \forall \,|x|\ge M, 
\end{equation}
{\em NB:} Of course, 
$$
R_+ \le r_+, \quad R_- \ge  r_-.
$$
This condition is the natural part of (\ref{b2}) in a combination with (\ref{b}); so, it will not be repeated each time in the statements. 

\medskip

The next theorem is the main result of the paper.
\begin{theorem}\label{thm2}
Let conditions (\ref{al}), (\ref{b}), (\ref{c2a}), (\ref{b2})
be satisfied. 
Then the process $(X,Z)$ is 2-recurrent, that is, 
there exists $C>0$ such that for any $M_1$ large enough and all $x \in \mathbb R^d$ and for $z=0,1$
\begin{equation}\label{e3a}
\mathsf E_{x,z}\tau^2_{M_1} 
\le C (x^{6} + 1). 
\end{equation} 
\end{theorem}

\ifhide
\medskip 

Under slightly more restrictive conditions it is possible to establish a bit better bound in comparison to (\ref{e3a}). Namely, assume 

\begin{equation}\label{c2R}
x b_-(x) \ge - R_-, \quad x b_+(x)\ge + R_+, \quad \forall \,|x|\ge M,
\end{equation}
along with
\begin{equation}\label{c2R2}
R_+\le r_+, \quad R_-\ge r_-.
\end{equation}
The following version of theorem 1 holds. 

\begin{theorem}\label{thm3}
Let conditions (\ref{al}), (\ref{b}), (\ref{c2}), (\ref{c2R}), (\ref{c2R2})
be satisfied. 
Then the process $(X,Z)$ is 2-recurrent, that is, 
there exists $C$ such that for any $M_1$ large enough and all $x \in \mathbb R^d$ and for $z=0,1$
\begin{equation}\label{et2-1}
\mathsf E_{x,z}\tau^2_{M_1} 
\le C (x^{6} + 1). 
\end{equation} 
\end{theorem}

\fi

\begin{Remark}\label{remthm2}
It is likely that conditions (\ref{al}), (\ref{b}), (\ref{c2}), (\ref{b2})  
suffice for the 2-recurrence of the process $(X,Z)$: namely, 
there exist $C, \delta>0$ such that for any $M_1$ large enough and all $x \in \mathbb R^d$ and for $z=0,1$
\begin{equation}\label{e3b}
\mathsf E_{x,z}\tau^2_{M_1} 
\le C (x^{4+\delta} + 1). 
\end{equation} 
(Of course, the less $\delta$ the greater $C$ here.) 
An accurate justification of this bound will be one of the goals, among others, for the further studies.
\end{Remark}

\begin{Remark}[On convergence]\label{remthm22}
The bound (\ref{e3b}) leads to the rate of convergence of the marginal distributions of the couple $(X_t,Z_t)$ in total variation towards the unique invariant measure via the coupling inequality:
\begin{equation}\label{mutinftybd}
\|\mu_t^{x,z} - \mu_\infty\|_{TV} \le \frac{C(x^6+1)}{1+t^2}, \quad t\ge 0.
\end{equation}
In principle, such a link is well-known in the recurrent Markov processes theory. The estimate (\ref{mutinftybd}) is particularly useful for the probabilistic representations of solutions of Poisson equations in the whole space $\mathbb R^d$. We postpone the details and the proof till futher publications.

\end{Remark}

\section{Auxiliaries I}\label{sec_auxI} %(for $\E \tau$)}
\label{S:2}
Here we state some auxiliary results similar to those established earlier in \cite{Ver21,Ver21b}, although, slightly different, which will be used in what follows. 

~

\noindent
Denote $\|b\| = \sup_{x,z}|b(x,z)|$.
Let $M_1 \gg M$ (the value $M_1$ to be specified later). 
Let
$$
T_0 : = \inf(t\ge 0: Z_t = 0), 
$$
and   
$$0 \le T_0 < T_1 < T_2 < \ldots, $$
where $T_n$ for each $n\ge 1$ is defined by induction as
 
$$
T_{n}:= \inf(t>T_{n-1}: \, Z_{t}-Z_{t-} \neq 0). 
$$
Let 
\begin{equation}\label{deftau}
\tau : = \inf(T_{2n}: \, |X_{T_{2n}}|\le M_1).
\end{equation}
{\em NB: In \cite{Ver21b} a bit different definition for $\tau$ was used, namely, $ \inf(T_{n}\ge 0: \, |X_{T_{n}}|\le M_1)$; in the present setting, given the new goal, the version (\ref{deftau}) is more convenient, as it is now not suitable to stop at ``odd'' stopping times $T_{2n+1}$.} 

To prove the main theorem (not in this section) it suffices to evaluate from above the value $\mathsf E_{x,z}\tau^2$ because apparently $\tau_{M_1} \le \tau$. 
%{\bf Assume $T_0=0$.}

Let $\epsilon>0, 0<q<1$ be any values satisfying the equality 
\begin{equation}\label{lle}
 \lambda_- (2r_++d +\epsilon) = q  \lambda_+ (2r_- -d - \epsilon)
\end{equation}
(see (\ref{c1})).
%{\color{green}
%\[
%\frac{(2r_- -  d)}{\lambda_-} > \frac{(2r_++ d)}{{\lambda}_+} \leqno{(\ref{c1})}
%\]
%}
Also, in the proof of the theorem it suffices to assume $|x| > M$ and  even $|x| > M_1$. The proofs of most of the lemmata and corollaries in this section follow from \cite{Ver21b}; if not, they are provided in what follows. 

\begin{lemma}\label{lem1}
1. Under the assumptions (\ref{al}) and (\ref{c1}) 
for any $\delta >0$ there exists $M_1$ such that 
\begin{align}\label{eps}
\max\left[\sup_{|x|>M_1}\mathsf E_{x,z} \left(\!\!\int_0^{T_1}1(\inf_{0\le s\le t}|X_s| \!\le\! M)dt|Z_0\!=\!0\!\right)\!, 
 \right. \nonumber\\ \\ \nonumber\left.
 \!  \sup_{|x|>M_1}\mathsf E_{x,z} \!\left(\!\int_0^{T_0}\!1(\inf_{0\le s\le t}|X_s| \!\le\! M)dt|Z_0\!=\!1\!\right)\right]
<\delta .
\end{align}
2. More than that, for any $\delta >0$ there exists $M_1$ such that 
\begin{align}\label{eps2}
\max\left[\sup_{|x|>M_1}\mathsf E_{x,z} \left(\!\!\int_0^{T_2}1(\inf_{0\le s\le t}|X_s| \!\le\! M)dt|Z_0\!=\!0\!\right)\!, 
 \right. \nonumber\\ \\ \nonumber\left.
 \!  \sup_{|x|>M_1}\mathsf E_{x,z} \!\left(\!\int_0^{T_1}\!1(\inf_{0\le s\le t}|X_s| \!\le\! M)dt|Z_0\!=\!1\!\right)\right]
<\delta .
\end{align}
\end{lemma}
{\em Proof.} 1. The proof of the first part of this lemma (i.e., of the bound (\ref{eps})) may be found in \cite{Ver21b}. 

\medskip

\noindent
2. The proof of (\ref{eps2}) is quite similar: the only change is that one random interval $[0,T_1]$ must be replaced by the union of two: 
$[0,T_2]=[0,T_1] \cup [T_1,T_2]$, which does not affect the calculus too much. Lemma \ref{lem1} follows. \hfill QED

\medskip

Let us denote by $X^i_t, \, i=0,1$ the solutions of the equations 
\begin{equation}\label{sde0}
dX^i_{t} =b(X^i_{t}, i)\, dt+ 
dW_{t}, \quad t\ge 0,
\quad X^i_{0} =x. 
\end{equation}
The index $z$ may be dropped for expectations and probabilities for each $X^i$.

\begin{lemma}\label{lem2}
If $M_1$ is large enough, then under the assumptions  (\ref{al}) and (\ref{c1}) 
for any $|x|>M_1$ for any $k=0,1,\ldots$
\begin{align}
\mathsf E_{x,z} (X_{T_{2k+1}\wedge   \tau}^2|{\mathcal F}_{T_{2k}}) 
\le \mathsf E_{x,z} (X_{T_{2k}\wedge   \tau}^2|{\mathcal F}_{T_{2k}}) 
 \nonumber \\  \nonumber %\label{ele2a} 
 \\
- 1(\tau > T_{2k}) \E(T_{2k+1}\wedge   \tau - T_{2k}\wedge   \tau |{\mathcal F}_{T_{2k}})((2r_--d)- \epsilon)
% \le x^2 
 \nonumber \\ \nonumber \\ \label{ele2a0}
\le \mathsf E_{x,z} (X_{T_{2k}\wedge   \tau}^2|{\mathcal F}_{T_{2k}})- 1(\tau > T_{2k}) \lambda_-^{-1}((2r_--d)- \epsilon); 
\end{align}
also,
\begin{align}\label{ele2b}
&\mathsf E_{x,z} (X_{T_{2k+2}\wedge   \tau}^2|{\mathcal F}_{T_{2k+1}})  
\le  \mathsf E_{x,z} (X_{T_{2k+1}\wedge   \tau}^2  |{\mathcal F}_{T_{2k+1}}) 
 \nonumber\\ \nonumber \\%\label{ele2b0}
&+  1(\tau > T_{2k}) \E(T_{2k+2}\wedge   \tau-T_{2k+1}\wedge   \tau | {\mathcal F}_{T_{2k+1}})((2r_{+}+d)+ \epsilon)
 \nonumber \\ \nonumber \\ %\label{ele2B0}
&\le \mathsf E_{x,z} (X_{T_{2k+1}\wedge   \tau}^2  |{\mathcal F}_{T_{2k+1}}) 
+ 1(\tau > T_{2k})  {\lambda}_+^{-1}((2r_++d)+ \epsilon).
\end{align}
\end{lemma}
Along with next lemma \ref{lem3}, this lemma \ref{lem2} will be used later on in lemma~\ref{lem11}. Note that with a new definition of $\tau$ the identity holds
\begin{equation}\label{tau2k}
1(\tau > T_{2k+1}) = 1(\tau > T_{2k}), 
\end{equation}
because $\tau$ may only be equal to some $T_{n}$ with an even index,  $n=2k$.

\medskip

\noindent
{\em Proof.} The proof of (\ref{ele2a0}) coincides verbatim with the proof of the bound (12) in \cite[Lemma 2]{Ver21b}, even though the setting in \cite{Ver21b} is a bit more general in the sense that the intensities may be variable with a possible dependence on the component $X_t$ at each time. The proof of (\ref{ele2b}) is also similar, but there is some difference because of the ``another indicator'' $1(\tau>T_{2k})$ in place of  $1(\tau>T_{2k+1})$. Let us show the full proof of (\ref{ele2b}). It suffices to consider the case $k=0$.

Under the condition $Z_0=1$ the process $X_t$ coincides with $X^1_t$ until the moment $T_0$. Hence, we have on $[0,T_0]$ as well as on $[T_1,T_2]$ by It\^o's formula 
\begin{align}\label{x2ito}
dX_t^2 - 2X_t dW_t = (2X_t b_+(X_t) + d)\,dt \le (2r_+ + d)dt,
\end{align}
the latter inequality on the set $(|X_t|> M)$ due to the assumptions (\ref{b}). Further, since $1(|X_t| > M) = 1 - 1(|X_t| \le M)$, we obtain
\begin{align*}
&\int_{T_1\wedge   \tau}^{T_2\wedge   \tau} 2X_t b_+(X_t)dt  
 \\\\
&= 
\int_{T_1\wedge   \tau}^{T_2\wedge   \tau}  2X_t b_+(X_t) 1(|X_t| > M)dt  
+\int_{T_1\wedge   \tau}^{T_2\wedge   \tau}  2X_t b_+(X_t)1(|X_t| \le M)dt 
 \\\\
&\le 2r_+ \int_{T_1\wedge   \tau}^{T_2\wedge   \tau} 1(|X_t| > M)dt  
+\int_{T_1\wedge   \tau}^{T_2\wedge   \tau}  2M \|b\| 1(|X_t| \le M)dt 
 \\\\ 
&= 2r_+ \int_{T_1\wedge   \tau}^{T_2\wedge   \tau} 1dt  
+\int_{T_1\wedge   \tau}^{T_2\wedge   \tau}  (2M \|b\| - 2r_+) 1(|X_t| \le M)dt 
 \\\\ 
&\le  2r_+ \int_{T_1\wedge   \tau}^{T_2\wedge   \tau} 1dt  
+2M \|b\|\int_{T_1\wedge   \tau}^{T_2\wedge   \tau}  1(|X_t| \le M)dt.
\end{align*}
Thus,  for $|x|>M_1$ (which is equivalent to $\tau > T_0$ in the case $T_0=0$) we get,
\begin{align*}
&\mathsf E_{x,z} \int_{T_1\wedge   \tau}^{T_2\wedge   \tau}  2X_t b_+(X_t)dt  = \mathsf E_{x} \int_{T_1\wedge   \tau}^{T_2\wedge   \tau}  2X^1_t b_+(X^1_t)dt 
 \\\\
&\le 2r_+ \E_x\int_{T_1\wedge   \tau}^{T_2\wedge   \tau} 1dt
+ 2M \|b\|  \E_x\int_{T_1\wedge   \tau}^{T_2\wedge   \tau}  1(|X^1_t| \le M)dt 
\\\\
&= 2r_+ \mathsf E_x\int_{T_1\wedge   \tau}^{T_2\wedge   \tau} 1dt
+ 2M \|b\| \mathbb  E_x\int_{T_1\wedge   \tau}^{T_2\wedge   \tau} 1(|X^1_t| \le M)dt 
 \\\\
&\le 2r_+ \mathsf E_x\int_{T_1\wedge   \tau}^{T_2\wedge   \tau} 1dt 
+ 2M \|b\|  \mathsf E_x\int_{T_1\wedge   \tau}^{T_2\wedge   \tau}  1(|X^1_t| \le M)dt 
  \\\\
&\le 2 r_+ \mathsf E_x\int_{T_1\wedge   \tau}^{T_2\wedge   \tau} 1dt 
+ 2M \|b\|\delta .
\end{align*}
For a fixed $\epsilon>0$ let us choose $\delta  = \lambda_+^{-1}\epsilon / (2M \|b\|) $. Then we get due to (\ref{x2ito}), 
\begin{align*}
\mathsf E_{x,z} (X_{T_2\wedge   \tau}^2\vert {\mathcal F}_{T_{1}}) - X_{T_{1}}^2  
\le (2r_++d)\mathsf E_x \left(\int^{T_2}_{T_1} dt \vert {\mathcal F}_{T_{1}}\right) 
+ \lambda_+^{-1}\epsilon
= \lambda_+^{-1}((2r_++d)+ \epsilon).
\end{align*}
Substituting here $X_{T_{2k+1}}$ instead of $X_{T_{1}}$, and writing $\mathsf E_x(\cdot |{\mathcal F}_{T_{2k+1}})$ instead of $\mathsf E_x(\cdot |{\mathcal F}_{T_{1}})$, and multiplying by $1(\tau > T_{2k+1}) = 1(\tau > T_{2k})$ (see (\ref{tau2k})), we obtain the bound (\ref{ele2b}), as required. \hfill QED

\begin{Corollary}\label{Cor2}
If $M_1$ is large enough, then under the assumptions  (\ref{al}) and (\ref{c1})
for any $|x|>M_1$ for any $k=0,1,\ldots$
\begin{align*}
&\mathsf E_{x,0} X_{T_{2k+1}\wedge   \tau}^2 - 
\mathsf E_{x,0}  X_{T_{2k}\wedge   \tau}^2
 \\\\
&\le - \mathsf E_{x,0} 1(\tau > T_{2k}) \E_{x,0}(T_{2k+1}\wedge   \tau - T_{2k}\wedge   \tau |{\mathcal F}_{T_{2k}})((2r_--d)- \epsilon)
 \\\\ 
&= - \mathsf E_{x,0} (T_{2k+1}\wedge   \tau - T_{2k}\wedge   \tau)((2r_--d)- \epsilon)
 \\\\
&\le - \mathsf E_{x,0} 1(\tau > T_{2k}) \lambda_-^{-1}((2r_--d)- \epsilon),
\end{align*}
and
\begin{align*}
&\mathsf E_{x,1} X_{T_{2k+2}\wedge   \tau}^2 - \mathsf E_{x,1}  X_{T_{2k+1}\wedge   \tau}^2  
 \nonumber\\ \\
&\le  \E_{x,1} 1(\tau > T_{2k}) (T_{2k+2}\wedge   \tau-T_{2k+1}\wedge   \tau ) ((2r_{+}+d)+ \epsilon)
 \nonumber\\ \\
&=  \E_{x,1} (T_{2k+2}\wedge   \tau-T_{2k}\wedge   \tau ) ((2r_{+}+d)+ \epsilon)
 \nonumber \\ \nonumber \\ 
&\le \mathsf E_{x,1} 1(\tau > T_{2k})  {\lambda}_+^{-1}((2r_++d)+ \epsilon).
\end{align*}
\end{Corollary}

\begin{lemma}\label{lem3}
If $M_1$ is large enough, then under the assumptions  (\ref{al}) and (\ref{c1}) 
for any $k=0,1,\ldots$ 
\begin{align}\label{ele3}
&\mathsf E_{x,z} (X_{T_{2k+2}\wedge   \tau}^2 |{\mathcal F}_{T_{2k+1}}) 
\le  \mathsf E_{x,z} (X_{T_{2k+1}\wedge   \tau}^2  | {\mathcal F}_{T_{2k+1}})  
 \nonumber \\ \nonumber\\
&+  1(\tau > T_{2k}) \mathsf E_{x,z} (T_{2k+2}\wedge \tau - T_{2k+1}\wedge \tau   |{\mathcal F}_{T_{2k+1}}) ) 
((2r_++d) + \epsilon))
 \nonumber \\  \nonumber\\ 
&\le \mathsf E_{x,z} (X_{T_{2k+1}\wedge   \tau}^2  |{\mathcal F}_{T_{2k+1}})   +  1(\tau > T_{2k})  {\lambda}_+^{-1}((2r_++d) + \epsilon)), 
\end{align}
and 
\begin{align}\label{ele7}
&\mathsf E_{x,z} (X_{T_{2k+1}\wedge   \tau}^2 | {\mathcal F}_{T_{2k}}) 
\le  \mathsf E_{x,z} (X_{T_{2k}\wedge   \tau}^2  |{\mathcal F}_{T_{2k}})  
  \nonumber \\ \nonumber\\ \nonumber 
&+  1(\tau > T_{2k}) \mathsf E_{x,z} (T_{2k+1}\wedge \tau - T_{2k}\wedge \tau   |{\mathcal F}_{T_{2k}}) ) 
  \nonumber \\ \nonumber \\ 
&\le  \mathsf E_{x,z} (X_{T_{2k}\wedge   \tau}^2  |{\mathcal F}_{T_{2k}})  -  1(\tau > T_{2k}) {\lambda}_-^{-1}((2r_--d) - \epsilon)). 
\end{align}
Also, 
\begin{align}\label{ele333}
&\mathsf E_{x,z} (X_{T_{2k+2}\wedge   \tau}^2 | {\mathcal F}_{T_{2k}}) - \mathsf E_{x,z} (X_{T_{2k}\wedge   \tau}^2  | {\mathcal F}_{T_{2k}})
  \nonumber \\ \nonumber \\ 
&\le  +  1(\tau > T_{2k})  \left({\lambda}_+^{-1}((2r_++d) + \epsilon))  -  {\lambda}_-^{-1}((2r_--d) - \epsilon))\right).
\end{align}
\end{lemma}
Along with lemma \ref{lem2}, this lemma \ref{lem3} will be used in what follows in lemma \ref{lem11}.

\begin{Corollary}\label{Cor3}
If $M_1$ is large enough, then under the assumptions   (\ref{al}) and  (\ref{c1}) 
for any $k=0,1,\ldots$ 
\begin{align*}%\label{ele3}
&\mathsf E_{x,0} X_{T_{2k+2}\wedge   \tau}^2 - \mathsf E_{x,0} X_{T_{2k+1}\wedge   \tau}^2  
 \nonumber \\\\
&\le  \E_{x,0} 1(\tau > T_{2k}) (T_{2k+2}\wedge \tau - T_{2k+1}\wedge \tau) 
((2r_++d) + \epsilon)
 \nonumber \\\\
&=  \E_{x,0} (T_{2k+2}\wedge \tau - T_{2k+1}\wedge \tau) 
((2r_++d) + \epsilon)
 \nonumber \\  \nonumber\\ %\label{ele30} 
&\le \mathsf E_{x,0}  1(\tau > T_{2k})  {\lambda}_+^{-1}((2r_++d) + \epsilon)), 
\end{align*}
and 
\begin{align*}%\label{ele7}
&\mathsf E_{x,1} X_{T_{2k+1}\wedge   \tau}^2 - 
\mathsf E_{x,1} X_{T_{2k}\wedge   \tau}^2  
  \nonumber \\\\ \nonumber 
&\le \mathsf E_{x,1} 1(\tau > T_{2k}) (T_{2k+1}\wedge \tau - T_{2k}\wedge \tau) 
  \nonumber \\ \nonumber \\ %\label{ele70}
&\le   -  \mathsf E_{x,1}1(\tau > T_{2k}) {\lambda}_-^{-1}((2r_--d) - \epsilon)). 
%\label{ele7}
\end{align*}
\end{Corollary}

\begin{lemma}\label{lem4}
Under the assumptions   (\ref{al}) and  (\ref{c1}) 
for any $k=0,1,\ldots$
\begin{align*}
 1(\tau>T_{2k}) \mathsf E_{X_{T_{2k+1}},1} (T_{2k+2}\wedge \tau - T_{2k+1}\wedge \tau) \ge  1(\tau>T_{2k})    \lambda_+^{-1},
\end{align*}
and 
\begin{align*}
 1(\tau>T_{2k}) \mathsf E_{X_{T_{2k}},0} (T_{2k+1}\wedge \tau - T_{2k}\wedge \tau) \le  1(\tau>T_{2k})     \lambda_-^{-1},
\end{align*}
\end{lemma}

\begin{lemma}\label{lem44}
Under the assumptions   (\ref{al}) and  (\ref{c1}) 
for any $k=0,1,\ldots$

\begin{align}\label{lem44-1}
&\mathsf E_{x,0} X_{T_{2k+2}\wedge \tau}^2 - \mathsf E_{x,0} X_{T_{2k}\wedge   \tau}^2  
\equiv \mathsf E_{x,0} 1(\tau > T_{2k}) (X_{T_{2k+2}\wedge \tau}^2 - \mathsf E_{x,0} X_{T_{2k}\wedge \tau}^2)
 \nonumber \\\\ \nonumber
&\le \mathsf E_{x,0}  1(\tau > T_{2k})  ({\lambda}_+^{-1}((2r_++d) + \epsilon)) - 
{\lambda}_-^{-1}((2r_--d) + \epsilon))
\end{align}
\end{lemma}
{\em NB:} The lemmata in this section suffice for an upper bound of $\E\tau$. They will also be used in the second order bound in the next sections.

\ifproofs

\noindent
{\em Proof of lemma \ref{lem4}.}
On the set $\tau>T_{2k+1}$ we have, 
\begin{align*}
 \mathsf E_{X_{T_{2k+1}},1} (T_{2k+2}\wedge \tau - T_{2k+1}\wedge \tau) 
= \mathsf E_{X_{T_{2k+1}},1} (T_{2k+2} - T_{2k+1}) 
\in [  \lambda_1^{-1} ,  \lambda_1^{-1}].
\end{align*}
Similarly, on the set $\tau>T_{2k}$
\begin{align*}
 \mathsf E_{X_{T_{2k}},0} (T_{2k+1}\wedge \tau - T_{2k}\wedge \tau) 
= \mathsf E_{X_{T_{2k}},0} (T_{2k+1} - T_{2k}) 
\in [  \lambda_-^{-1} ,  \lambda_-^{-1}].
\end{align*}
On the sets $\tau\le T_{2k+1}$ and $\tau\le T_{2k}$, respectively, both sides of the required inequalities equal zero.  Lemma \ref{lem4} follows. \hfill QED

\fi

\section{Auxiliaries II: fourth moments of $X_t$}\label{sec_auxII} 
In this section new auxiliary results will be established which will be required for the proof of the main theorem. 

~

\noindent
Potentially, the recurrence in terms of the $4+\delta$ moment of the  marginals of the process  $X$  should lead to a bound for the second order moment of $\tau$; the author is not sure whether the exact forth moment suffices, although, it might be. This section is devoted to the bounds for the fourth moment of $X$, and the sixth moment will be esstimated in the next section. Recall that the additional stronger assumptions (\ref{c2}), (\ref{b2})
are required for this section. 
\ifhide
{\color{green}
\[
\frac{(4r_--6d)}{\lambda_-} > \frac{(4r_++6d)}{\lambda_1},
\leqno{(\ref{c2})}
\]
\[
x b_-(x) \ge - R_-, \quad x b_+(x)\ge + R_+, \quad \forall \,|x|\ge M.
\leqno{(\ref{b2})}
\]
}
\fi
We are going to evaluate from below and from above the differences  of the fourth moment of the diffusion component $X_t$ on ``positive'' and ``negative'' intervals between the consequtive jumps of the $Z$-component, i.e., where $Z=0$ and $Z=1$.

\begin{lemma}\label{lem50}
Let the assumptions (\ref{b}) and (\ref{b2}) be satisfied. Then  there exists $C>0$ such that for any $x,z,t$
\begin{align}\label{l50-1}
&\E_{x,z} X_t^2 \le x^2 + Ct, 
 \\ \nonumber \\ \label{l50-2}
&\E_{x,z} X_t^2 \ge x^2 - Ct,  
\\  \nonumber \\ \label{l50-3}
&\E_{x,z} X_t^4 
\le  x^4 + Cx^2t + Ct^2 +Ct,
 \\  \nonumber \\ \label{l50-4}
&\E_{x,z} X_t^4 
\ge x^4  - Cx^2t - Ct^2 -Ct.
\end{align}

\end{lemma}
{\em NB:} For what follows it is important/useful to keep the senior coefficients equal to one in all right hand sides here: $1\times x^2$ in (\ref{l50-1}), $1\times x^4$ in (\ref{l50-3}), et al. The constants $C$ in the right hand sides are less important for what follows.

~

\noindent
{\em Proof.} 
\ifhide
Firstly notice that
\begin{align}\label{l50proof-1}
\E_{x,z} |X_t| \le |x|+ Ct + d\sqrt{t}.
\end{align}
(NB: In the latter formula $d$ is the dimension.) Then, 
\fi
By It\^o's formula, 
\begin{equation}\label{itox2}
dX^2_t = 2X_tdX_t + (dX_t)^2 = (2X_tb(X_t,Z_t)+d)dt  + 2X_t dW_t,
\end{equation}
and 
\begin{align*}%\label{itox4}
&dX^4_t =  2X^2_t [(2X_t b(X_t,Z_t) + d + 2)\,dt + 2X_t dW_t].
\end{align*}
%and for the sequel let us note that 
%\begin{align*}%\label{itox4}
%&dX^6_t =  3X^4_t [(2X_t b(X_t,Z_t) + d + 4)\,dt + 2X_t dW_t] 
%\end{align*}

\noindent
Note that due to the assumptions (\ref{b}) and (\ref{b2}), and by virtue of the local boundedness of the drift we have 
\begin{equation}\label{l50proof-3}
-\infty < \inf_{x,z}xb(x,z) \le 
\sup_{x,z}xb(x,z)< \infty.
\end{equation}
So, 
\begin{align*}
\E_{x,z} X_t^2 \!-\! x^2 
\!=\! 2\E_{x,z} \int_0^t (X_s b(X_s,Z_s) \!+\! d)ds 
\!=\! 2 \int_0^t \E_{x,z} \underbrace{(X_s b(X_s,Z_s)\!+\!d)}_{\le C}ds \!\le\! Ct, 
\end{align*}
which shows (\ref{l50-1}).

To justify (\ref{l50-2}) we conclude using (\ref{itox2}) 
\begin{align*}
&\E_{x,z} X_t^2\! -\! x^2 \!=\! 2\E_{x,z} \int_0^t (X_s b(X_s,Z_s)\!+\!d)ds \!=\! 2 \int_0^t \E_{x,z} \underbrace{(X_s b(X_s,Z_s)\!+\!d)}_{\ge -C}ds \!\ge \!-\!Ct. 
\end{align*}
\ifhide
\begin{align*}
&\E_{x,z} X_t^2 - x^2 = 2\E_{x,z} \int_0^t X_s b(X_s,Z_s)ds = 2 \int_0^t \E_{x,z} X_s b(X_s,Z_s)ds 
 \\\\
&\ge - C\int_0^t \E_{x,z}|X_s|ds \ge -C\int_0^t (|x|+ Cs + C\sqrt{s})ds
 \\\\
&= -C|x|t - Ct^2 - Ct^{3/2} \ge -C|x|t - Ct^2 - Ct^{} .
\end{align*}
\fi
It is needless to say that all constants $C$ are or may be different here. %The change from $t^{3/2}$ is a bit artificial and not important for what follows: just further calculi will look a bit less involved.
Now (\ref{itox2}) and (\ref{l50-1}) imply (\ref{l50-3}) as follows, 
\begin{align*}
&\E_{x,z} X_t^4 
\le x^4 + C\E \int_0^t (X_s^2X_sb(X_s,Z_s) +C)ds
\le x^4 + C \int_0^t (\E X_s^2 +C)ds 
 \\\\
&\le x^4 + C \int_0^t (x^2 + Cs +C)ds 
\le  x^4 + Cx^2t + Ct^2 +Ct,
\end{align*}
as required.
Similarly, 
\begin{align*}
&\E_{x,z} X_t^4 
\ge x^4 + 4\E \int_0^t (X_s^2X_sb(X_s,Z_s) - C)ds
\stackrel{(\ref{l50proof-3})}\ge x^4 - C \int_0^t (\E X_s^2 +C)ds 
 \\\\
&\ge x^4 - C \int_0^t (x^2 + Cs +C)ds 
\ge  x^4 - Cx^2t - Ct^2 -Ct,
\end{align*}
which shows (\ref{l50-4}). Lemma \ref{lem50} is proved. 
\ifhide
Let us show (\ref{l50-4}). Using 
$$
\E_{x,z} |X|_t^3\le (\E_{x,z} X_t^4)^{3/4} \le 
|x|^3 + C|x|^{3/2} t + Ct^{3/2} +Ct^{3/4} \le 
$$
due to (\ref{l50-3}), we estimate, 
\begin{align*}
&\E_{x,z} X_t^4 \ge x^4 - C\E \int_0^t (X_s^3b(X_s,Z_s) +C)ds
 \\\\
&\ge x^4 - C\int_0^t (\E |X_s^3|  +C)ds
 \\\\
&\ge x^4 - C\int_0^t (|x|^3 + C|x|^{3/2} s + Cs^{3/2} +Cs^{3/4}+C)ds
 \\\\ 
& = x^4 - C t |x|^3 - C|x|^{3/2} t^2 - Ct^{5/2} -Ct^{7/4}-C t
\\\\ 
& \ge  x^4 - C t |x|^3 - C|x|^{3/2} t^2
- Ct^{3} -Ct^{2}-C t, 
\end{align*}
as required??
\fi
\hfill QED

\begin{lemma}\label{lem5}
%(On + interval) 
If $|x|>M_1$, then there exists $C>0$ such that on the positive interval under the assumptions (\ref{al}),  (\ref{b}),  (\ref{c2}), and (\ref{b2})
\begin{align}\label{x4+}
\E_{x,1} (X_T^4 - x^4) \le
\frac{(4r_++(2d+4)) }{\lambda_+} x^2 + C,
\end{align}
or, equivalently, 
\begin{align}\label{lem7x44+}
\E_{x,1} (X_{T\wedge\tau})^4 - x^4 \le
\frac{(4r_++(2d+4)) }{\lambda_+} x^2 + C.
\end{align}

\end{lemma}
\noindent
{\em Proof.} For $|x|>M_1$ we have
\begin{align*}
&\E_{x,1} (X_T^4 - x^4) = \E_{x,1}\int_0^T (4 \underbrace{X_s^3b(X_s)}_{\le r_+ X_s^2 1(X_s^2>M^2)+C, \, \ge R_+X_s^2 1(X_s^2>M^2)-C} + (2d+4)X_s^2)ds 
 \\\\
&\le (4r_++6) \E_{x,1}\int_0^T  (X_s^2 +C)ds 
= (4r_++(2d+4)) \E_{x,1}\int_0^T  (X_{s\wedge\tau}^2 +C)ds 
 \\\\
&= (4r_++(2d+4)) \E_{x,1}\int_0^\infty 1(s<T)  (X_{s\wedge\tau}^2+C) ds 
 \\\\
&= (4r_++(2d+4)) \int_0^\infty \E_{x,1} 1(s<T)  (X_s^2+C) ds 
 \\\\
&= (4r_++(2d+4)) \int_0^\infty E_x 1(s<T)  \E_{x,1} (X_{s}^2 +C)ds 
 \\\\
&\le (4r_++(2d+4)) \int_0^\infty \exp(-\lambda_+ s)  (x^2+ Cs) ds 
+\frac{C}{\lambda_+}
 \\\\
&= \underbrace{\frac{(4r_++(2d+4)) }{\lambda_+}}_{=:\kappa_+} x^2 + C=: \kappa_+ x^2 + C,
\end{align*}
with some $C$. 
This justifies the estimate (\ref{x4+}). The lemma is proved. \hfill QED

%\newpage

\begin{lemma}\label{lem8}
If $|x|>M_1$,  then there exists $C>0$ such that on the negative interval under the assumptions (\ref{al}),  (\ref{b}),  (\ref{c2}), and (\ref{b2})\begin{align}\label{l81}
\!\E_{x,0} (X_{T}^4 \!-\! x^4) \!=\!
\E_{x,0} (X_{T\wedge\tau}^4 \!-\! x^4) 
\!\le \! - \frac{(4r_-\!-\!(2d+4))}{\lambda_-} x^2 \!+ C. 
\end{align}
\end{lemma}
\noindent
{\em Proof.}
We have on the set $\tau>0$ %(on a negative interval) 
with $|x|>M_1$:
\begin{align*}
&\E_{x,0} (X_T^4 - x^4) = \E_{x,0}\int_0^T (4 \underbrace{X_s^3b(X_s)}_{\le - r_- X_s^2 +C 1(X_s^2\le M)} + (2d+4)X_s^2 )ds 
 \\\\
&\le (-4r_-+(2d+4)) \E_{x,0}\int_0^T  X_s^2 ds 
+ C \E_{x,0}\int_0^T  X_s^2 1(X_s^2\le M)ds + C
 \\\\
&= (-4r_-+(2d+4)) \E_{x,0}\int_0^\infty 1(s<T)  X_s^2 ds
 \\\\
&+ C \int_0^\infty \E_{x,0} 1(s<T) \E_x (X^0_s)^2 1((X^0_s)^2\le M^2)ds + C 
 \\\\
&= (-4r_-+(2d+4)) \int_0^\infty \E_{x,0} 1(s<T)  X_s^2 ds
 \\\\
&+ C \int_0^\infty \E_{x,0} 1(s<T) \E_{x,0} X_s^2 1(X_s^2\le M^2)ds + C
 \\\\
&= (-4r_-+(2d+4)) \int_0^\infty \E_{x,0} 1(s<T)  \underbrace{\E_{x,0} X_{s}^2}_{\ge x^2 + (-2R_- +1 -C)s} ds 
 \\\\
&+ C \underbrace{\int_0^\infty \exp(-\lambda_- s) \E_{x,0} X_s^2 1(X_s^2\le M^2)ds}_{\le M^2/\lambda_-} + C
 \\\\
&\le (-4r_-+(2d+4)) \int_0^\infty \exp(-\lambda_- s)  (x^2+(-2R_-+1 -C)s) ds 
 \\\\ 
&+ C \int_0^\infty \exp(-\lambda_- s) \E_{x,0} X_s^2 1(X_s^2\le M^2)ds + C
 \\\\
&= - \underbrace{\frac{(4r_--(2d+4))}{\lambda_-}}_{=:\kappa_-} x^2 + \frac{(4r_--(2d+4))(2R_--1+C)}{\lambda_-^2} +C
=: - \kappa_- x^2 + C,  
\end{align*}
with a certain $C$. This shows (\ref{l81}). 
Lemma \ref{lem8} is proved.  \hfill QED

\begin{lemma}\label{lem9}
Under the assumptions (\ref{al}), (\ref{b}), (\ref{c2}), and (\ref{b2}), it is possible to choose $M_1>\!\!>M$ so that 
$$\E_{x,0} X_T^4 \equiv \E_{x,0} X_{T_1}^4 \le x^4,$$ 
and, moreover, 
$$
\E_{x,0} X_{T_2\wedge \tau}^4 \le x^4 - cx^2.
$$ 
Similarly, 
$$
\E_{ X_{T_{2n}\wedge \tau}} 1(T_{2n}<\tau) X_{T_{2n+2}\wedge \tau}^4 \le 
\E_{X_{T_{2n}\wedge \tau}} 1(T_{2n}<\tau) X_{T_{2n}\wedge \tau}^4 - c\E_{X_{T_{2n}\wedge \tau}} 1(T_{2n}<\tau) X_{T_{2n}\wedge \tau}^2,
$$ 
and
$$
\E_{x,0} 1(T_{2n}<\tau) X_{T_{2n+2}\wedge \tau}^4 \le 
\E_x 1(T_{2n}<\tau) X_{T_{2n}\wedge \tau}^4 - c\E_x 1(T_{2n}<\tau) X_{T_{2n}\wedge \tau}^2,
$$ 
and, hence, by induction, 
$$
\E_{x,0} 1(T_{2n}<\tau) X_{T_{2n+2}\wedge \tau}^4 \le x^4 - cx^2.
$$ 
\end{lemma}
\noindent
{\em Proof.}
We evaluate after two steps (if $M_1$ is large enough)
\begin{align*}
&\E_{x,0} X_{T_2\wedge\tau}^4 - \E_{x,0}  X_{T_1\wedge\tau}^4 
= \E_{x,0}  \left(X_{T_2\wedge\tau}^4 -  X_{T_1\wedge\tau}^4\right) 
1(\tau > T_1) 
 \\\\
& = \E_{x,0}  1(\tau > T_1) \E_{x,0} \left(X_{T_2\wedge\tau}^4 -  X_{T_1\wedge\tau}^4 |{\cal F}_{T_1}\right)  
 \\\\
&\le \E_{x,0} 1(\tau > T_1)\left( \frac{(4r_++(2d+4)) }{\lambda_+}  X_{T_1\wedge\tau}^2 + \frac{(4r_++(2d+4)) (2r_++d+C)}{\lambda_+^2}  \right)
 \\\\
&\le \E_{x,0} 1(\tau > T_0)\left( \frac{(4r_++(2d+4)) }{\lambda_+}  X_{T_1\wedge\tau}^2 + \frac{(4r_++(2d+4)) (2r_++d+C)}{\lambda_+^2}  \right)
\\\\
&= \E_{x,0} 1(\tau > T_0)\frac{(4r_++(2d+4)) }{\lambda_+}  X_{T_1\wedge\tau}^2 + \frac{(4r_++(2d+4)) (2r_++d+C)}{\lambda_+^2}
 \\\\
& \le  \E_{x,0} 1(\tau > T_0)\frac{(4r_++(2d+4)) }{\lambda_+} \,x^2 
 +\E_{x,0} 1(\tau > T_0) {\lambda}_+^{-1}((2r_--d) - \epsilon))  
  \\\\
& + \frac{(4r_++(2d+4)) (2r_++d+C)}{\lambda_+^2},
\end{align*}
and
\begin{align*}
&\E_{x,0}  X_{T_1\wedge\tau}^4 - \E_{x,0}  X_{T_0\wedge\tau}^4 
= \E_x \left(X_{T_1\wedge\tau}^4 -  X_{T_0\wedge\tau}^4\right) 
1(\tau > T_0) 
 \\\\
&\le - \frac{(4r_--(2d+4))}{\lambda_-} x^2  
+C.
\end{align*}
Adding up, we get for $x^2>M_1^2$, if $M_1$ is large enough, 
\begin{align*}
\E_{x,0}  X_{T_2\wedge\tau}^4 - x^4
\le x^2 \left(- \frac{(4r_--(2d+4))}{\lambda_-} + \frac{(4r_++(2d+4))}{\lambda_+} \right) +C . 
\end{align*}
Due to the assumption (\ref{b2}), it is possible to choose $M_1$ so that 
$$
\frac12\,M_1^2 \left(\frac{(4r_--(2d+4))}{\lambda_-} - \frac{(4r_++(2d+4))}{\lambda_+}\right)  > C,
$$
with the value of the constant $C$ from the previous inequality. 
Then on $x^2>M_1^2$
\begin{align*}
\E_{x,0}  X_{T_2}^4 
\le  x^4 - cx^2,
\end{align*}
with a certain $c>0$, as required. The lemma follows. \hfill QED

\begin{Corollary}\label{cor33}
Under the assumptions assumptions (\ref{al}),  (\ref{b}),  (\ref{c2}), and (\ref{b2}), if $M_1$ is large enough, then  for any $k=0,1,\ldots$ 
\begin{align}\label{ele7a}
\mathsf E_{x,0} 1(T_{2k}<\tau) X_{T_{2k+2}\wedge   \tau}^4 
\le \mathsf E_{x,1} 1(T_{2k-2}<\tau) X_{T_{2k}\wedge   \tau}^4...
\le x^4 + C.
\end{align}

\end{Corollary}
{\em NB:} Similar inequalities may be also proved for the case $z=1$. However, since our ultimate goal is just an estimate for the second moment of $\tau$ (for any value of $z$), this is not necessary, as in the end the bound for $\E_{x,1}\tau^2$ will be derived from the bound on $\E_{x,0}\tau^2$.

\section{Auxiliaries III: sixth moments for $X_t$}\label{sec_auxIII}
The following is the analogue of lemma \ref{lem8} for $\E X_T^6$.
\begin{lemma}\label{lem8fr}
Suppose the assumptions (\ref{al}),  (\ref{b}),  (\ref{c2a}), and (\ref{b2}) hold true. In this case, 
if $|x|>M_1$,  then on the negative interval for some $C,C_->0$, 
\begin{align}\label{l10-1}
\!\E_{x,0} (X_{T}^{6} \!-\! x^{6}) \!=\! 
\!\E_{x,0} (X_{T\wedge\tau}^{6} \!-\! x^{6})
\!\le \! - \frac{(6r_-\!-\!(3d+12))}{\lambda_-} x^{4} \!+ C_-x^2 + C.
\end{align}
\end{lemma}

\noindent
{\em Proof.}
By It\^o's formula, 
\begin{align*}%\label{itox4}
&dX^6_t 
=  3X^4_t [(2X_t b(X_t,Z_t) + d + 4)\,dt + 2X_t dW_t]. 
\end{align*}
So, we have on the set $\tau>0$ %(on a negative interval) 
with $|x|>M_1$: 
\begin{align*}
&\E_{x,0} (X_T^6 - x^6) = \E_{x,0}\int_0^T (6 X_s^{4}\underbrace{X_s^{}b(X_s)}_{\le - r_- X_s^{4} +C 1(X_s^2\le M)} + (3d+12) X_s^{4})ds 
 \\\\
&\le (-6r_-+(3d+12) ) \E_{x,0}\int_0^T  X_s^{4} ds 
+ C \E_{x,0}\int_0^T  X_s^{4} 1(|X_s|\le M)ds
 \\\\
&\!=\! (\!-\!6r_-\!+\!(3d\!+\!12)) \E_{x,0}\int_0^\infty \!1(s\!<\!T)  X_s^{4} ds
\!+\! C \int_0^\infty \!\E_{x,0} 1(s\!<\!T) \E_x (X^1_s)^{4} 1((X^1_s)^2\!\le\! M^2)ds
 \\\\
&=\! (-6r_-\!+\!(3d\!+\!12) ) \int_0^\infty \E_{x,0} 1(s\!<\!T)  X_s^{4} ds
\!+\! C \int_0^\infty \E_{x,0} 1(s\!<\!T) \E_{x,0} X_s^{4} 1(X_s^2\!\le\! M^2)ds
 \\\\
&= (-6r_-+(3d+12) )  \int_0^\infty \E_{x,0} 1(s<T)  \underbrace{\E_{x,0} X_{s}^{4}}_{\ge x^4 - Cx^2 s -C s^2 - Cs} ds 
 \\\\
&+ C \underbrace{\int_0^\infty \exp(-\lambda_- s) \E_{x,0} X_s^{4} 1(X_s^2\le M^2)ds}_{\le M^{4}/\lambda_-}
 \\\\
&\le (-6r_-+(3d+12) )  \int_0^\infty \exp(-\lambda_- s)  (x^{4}
- Cx^2s - Cs^2 - Cs) ds 
+C
 \\\\ 
&+ C \int_0^\infty \exp(-\lambda_- s) \E_{x,0} X_s^{4} 1(X_s^2\le M^2)ds
 \\\\
&= - \underbrace{\frac{(6r_--(3d+12))}{\lambda_-}}_{=:\kappa_-} x^{4} 
+ C_-x^2 + C
=: - \kappa_- x^{4} + C_-x^2 + C.   
\end{align*}
This shows (\ref{l10-1}). 
Lemma \ref{lem8fr} is proved.  \hfill QED

\medskip

\ifp
{\color{green}
\begin{align*}%\label{l50-3}
&\E_{x,z} X_t^4 
\le  x^4 + Cx^2t + Ct^2 +Ct,
 \\  \nonumber \\ %\label{l50-4}
&\E_{x,z} X_t^4 
\ge x^4 
 - Cx^2t - Ct^2 -Ct,
\end{align*}
}
\fi

The next lemma is the analogue of lemma \ref{lem5}.

\begin{lemma}\label{lem5a}%lem11
If $|x|>M_1$, then on the positive interval under the assumptions (\ref{al}),  (\ref{b}),  (\ref{c2a}), and (\ref{b2}), for some $C,C_+>0$, 
\begin{align}\label{l11-1}
\E_{x,1} (X_T^6 - x^6) \le
\frac{(6r_++(3d+12)) }{\lambda_+} x^4 
+ C_+x^2 + C.
\end{align}

\end{lemma}
\noindent
{\em Proof.} For $|x|>M_1$ we have
\begin{align*}
&\E_{x,1} (X_T^6 - x^6) = \E_{x,1}\int_0^T (6 \underbrace{X_s^5 b(X_s)}_{\le r_+ X_s^4 1(X_s^2>M^2)+C, \, \ge R_+X_s^4 1(X_s^2>M^2)-C} + (3d+12) X_s^4)ds 
 \\\\
&\le (6r_++(3d+12) ) \E_{x,1}\int_0^T  (X_s^4 +C)ds 
= (6r_++(3d+12) ) \E_{x,1}\int_0^T  (X_{s\wedge\tau}^4 +C)ds 
 \\\\
&=\! (6r_+\!+\!(3d+12) ) \E_{x,1}\int_0^\infty \!1(s\!<\!T)  (X_{s\wedge\tau}^4\!+\!C) ds 
 \\\\ & 
\! = \!(6r_+\!+\!(3d+12) ) \int_0^\infty \!\E_{x,1} 1(s\!<\!T)  (X_s^4\!+\!C) ds 
 \\\\
&= (6r_++(3d+12) ) \int_0^\infty E_x 1(s<T)  \E_{x,1} (X_{s}^4 +C)ds 
 \\\\
&\stackrel{(\ref{l50-3})}\le (6r_++(3d+12) ) \int_0^\infty \exp(-\lambda_+ s)  (x^4 + Cx^2s + Cs^2 +Cs) ds 
+\frac{C}{\lambda_+}
 \\\\
&= \underbrace{\frac{(6r_++(3d+12) ) }{\lambda_+}}_{=:\kappa_+} x^4 + C=: \kappa_+ x^4 + C_+x^2 + C.
\end{align*}
This justifies the estimate (\ref{l11-1}). Lemma \ref{lem5a} follows. 
 \hfill QED

\begin{lemma}\label{lem9a}
Suppose that the assumptions (\ref{al}),  (\ref{b}),  (\ref{c2a}), and (\ref{b2}) hold true. Then it is possible to choose $M_1>\!\!>M$ so that 
$$
\E_{x,0} X_T^6 \equiv \E_{x,0} X_{T_1}^6 \le x^6,
$$ 
and, moreover, 
$$
\E_{x,0} X_{T_2\wedge \tau}^6 \le x^6 - cx^{4}.
$$ 
Similarly, 
$$
\E_{ X_{T_{2n}\wedge \tau}} 1(T_{2n}<\tau) X_{T_{2n+2}\wedge \tau}^6 \le 
\E_{X_{T_{2n}\wedge \tau}} 1(T_{2n}<\tau) X_{T_{2n}\wedge \tau}^6 - c\E_{X_{T_{2n}\wedge \tau}} 1(T_{2n}<\tau) X_{T_{2n}\wedge \tau}^{4},
$$ 
and
$$
\E_{x,0} 1(T_{2n}<\tau) X_{T_{2n+2}\wedge \tau}^6 \le 
\E_x 1(T_{2n}<\tau) X_{T_{2n}\wedge \tau}^6 - c\E_x 1(T_{2n}<\tau) X_{T_{2n}\wedge \tau}^{4},
$$ 
and (by induction) 
$$
\E_{x,0} 1(T_{2n}<\tau) X_{T_{2n+2}\wedge \tau}^6 \le x^6 - cx^{4}.
$$ 
\end{lemma}
\noindent
{\em Proof.}
We evaluate after two steps (if $M_1$ is large enough) on the positive interval
\begin{align*}
&\E_{x,0} X_{T_2\wedge\tau}^6 - \E_{x,0}  X_{T_1\wedge\tau}^6 
= \E_{x,0}  \left(X_{T_2\wedge\tau}^6 -  X_{T_1\wedge\tau}^6\right) 
1(\tau > T_1) 
 \\\\
& = \E_{x,0}  1(\tau > T_1) \E_{x,0} \left(X_{T_2\wedge\tau}^6 -  X_{T_1\wedge\tau}^6 |{\cal F}_{T_1}\right)  
 \\\\
&\le \E_{x,0} 1(\tau > T_1)\left( \frac{(6r_++(3d+12)) }{\lambda_+}  X_{T_1\wedge\tau}^4 + C X_{T_1\wedge\tau}^2 + C  \right)
 \\\\
&\le \E_{x,0} 1(\tau > T_0)\left( \frac{(6r_++(3d+12)) }{\lambda_+}  X_{T_1\wedge\tau}^4 + C X_{T_1\wedge\tau}^2 + C \right)
\\\\
&= \frac{(6r_++(3d+12)) }{\lambda_+}  \E_{x,0} 1(\tau > T_0)X_{T_1\wedge\tau}^4 + C\E_{x,0} 1(\tau > T_0)X_{T_1\wedge\tau}^2 +C
 \\\\
&\stackrel{\text{lemma \ref{lem5}}} \le  \frac{(6r_++(3d+12)) }{\lambda_+}  \E_{x,0} 1(\tau > T_0)\,(x^4 
 +  {\lambda}_+^{-1}((6r_++1) + \epsilon))x^2 + C)  
  \\\\
& \stackrel{\text{lemma \ref{lem3}}}+ \E_{x,0} 1(\tau > T_0) (x^2+C)
\end{align*}
(it was used that $1(\tau>T_1)=1(\tau>T_0)$), and on the negative interval
\begin{align*}
&\E_{x,0}  X_{T_1\wedge\tau}^6 - \E_{x,0}  X_{T_0\wedge\tau}^6 
= E_x \left(X_{T_1\wedge\tau}^6 -  X_{T_0\wedge\tau}^6\right) 
1(\tau > T_0) 
 \\\\
&\stackrel{\text{lemma \ref{lem8fr}}}\le  \E_{x,0} 1(\tau > T_0)
(- \frac{(6r_--(3d+12))}{\lambda_-} x^4 +  Cx^2 + C).
\end{align*}
Adding up, we get for $x^2>M_1^2$, if $M_1$ is large enough 
\begin{align}\label{Cx4x2}
&\E_{x,0}  X_{T_2\wedge\tau}^6 - x^6
\le x^4 \left(- \frac{(6r_--(3d+12))}{\lambda_-} + \frac{(6r_++(3d+12))}{\lambda_+} \right) + Cx^2 +C . 
\end{align}
We can choose $M_1$ so that 
$$
\frac12\, M_1^2 \left(\frac{(6r_--(3d+12))}{\lambda_-} - \frac{(6r_++(3d+12))}{\lambda_+}\right)  > C,
$$
for $C$ is from the term $Cx^2$ in (\ref{Cx4x2}).
Then on $x^2>M_1^2$
\begin{align*}
\E_{x,0}  X_{T_2}^6 
\le  x^6 - cx^4,
\end{align*}
with a certain $c>0$, as required. Lemma \ref{lem9a} is proved. \hfill QED

\begin{Corollary}\label{Cor4}
Under the assumptions (\ref{al}),  (\ref{b}),  (\ref{c2a}), and (\ref{b2}) for any $m\le 3$ and for any $n=0,1,\ldots$, 
$$
\E_{x,0} 1(T_{2n}<\tau) X_{T_{2n+2}\wedge \tau}^{2m} \le x^{2m}.
$$

\end{Corollary}
{\em Proof.} The bound follows straightforwardly from lemma \ref{lem9a} and from  H\"older's inequality. \hfill QED

\section{Second moment for $\tau_{}$ via $\E X_t^6$}\label{sec_mainproof}
Let
$$
Y_n:= X_{T_{2n}}, \quad N \,(=\tau^Y):=\inf(n: |Y_n|\le M_1).
$$
It is expected that under the appropriate conditions, the same as for $(X_t)$, the Markov process $(Y_n)$ is recurrent and ergodic and that the moment properties of the stopping time $N$ will help establish similar properties for the stopping time $\tau$ which is our ultimate goal. 
Indeed, 
$$
\tau = T_{2N} = T_0 + \sum_{i=0}^{N-1} \underbrace{(T_{2i+2} - T_{2i})}_{=: \eta_i}, 
$$
and the process 
$$
S_n:=  \sum_{i=0}^{n-1} (\eta_{i} - \E\eta_{i})
$$ is a square integrable discrete time martingale with the compensator
$$
\langle S\rangle_n = 2\left(\frac1{\lambda_-^2} + \frac1{\lambda_-\lambda_+} + \frac1{\lambda_+^2}\right)n =: c_2n.
$$
Note that $\E\eta_{i} = \lambda_-^{-1} +  \lambda_+^{-1}$.  
Hence, 
%due to Doob's optional stopping theorem, 
the moments of the random variable $\tau$ may be evaluated on the basis of the moments of the random variable $N$; in particular, since 
$$
\tau = T_0 + c_1 N + S_N, 
$$
we have a bound
\begin{align}\label{etau2}
&\E \tau^2 = \E (T_0+c_1 N + S_N)^2  \le 3\E T_0^2 + 3 c_1^2 \E N^2 + 3\E \langle S\rangle_N 
 \nonumber \\ \\ \nonumber
&= C + 3 (c_1^2 \E N^2 + c_2 \E N), 
\end{align}
and, more generally, for any $m>0$
$$  
\E \tau^m = \E (T_0 + c_1 N + S_N)^m  \le 3^{m-1} \E T_0^m + 3^{m-1} (c_1^m \E N^m + c_m^{m/2} \E N^{m/2}),
$$
with an appropriate $c_m$.
The bound on $\E N$ is based on the following 
\begin{lemma}\label{lem11}
Under the assumptions of theorem \ref{thm2}, for any $z=0,1$, 
\begin{equation}\label{lem11-1}
1(n+1<N)\E_{Y_{n}} |Y_{(n+1)}|^2 \le
1(n<N)\E_{Y_{n}} |Y_{(n+1)}|^2 \le 1(n<N)|Y_{n}|^2 - c  1(n<N).
\end{equation}
Also, for any $z=0,1$,
\begin{equation}\label{eN}
\E_{x,z} N \le C(x^2 +1).
\end{equation}

\end{lemma}
{\em Proof.} 
Consider the case $z=0$; then $T_0=0$ and $Y_0 = x$. The first inequality in (\ref{lem11-1}) is trivial because $1(n+1<N) \le
1(n<N)$.
Further, due to the inequality (\ref{ele333}), and  
\ifp
{\color{green}
\begin{align*}
\mathsf E_{x,z} (X_{T_{2k+2}\wedge   \tau}^2 | {\mathcal F}_{T_{2k}}) - \mathsf E_{x,z} (X_{T_{2k}\wedge   \tau}^2  | {\mathcal F}_{T_{2k}})
  \nonumber \\ \nonumber \\ %\label{ele70}
\le  +  1(\tau > T_{2k})  \left({\lambda}_+^{-1}((2r_++1) + \epsilon))  -  {\lambda}_-^{-1}((2r_--1) - \epsilon))\right).
\end{align*}
}
\fi
since the process $X_{T_{2k}}=Y_k$ is markovian, then (\ref{ele333}) may be rewritten as 
\begin{align*}
&\mathsf E_{x,z} (Y_{(k+1)\wedge  N}^2 | {\mathcal F}^Y_{{k}}) - \mathsf E_{x,z} (Y_{k\wedge  N}^2  | {\mathcal F}^Y_{k})
  \nonumber \\ \nonumber \\ 
&\le  +  1(N>k)  \left({\lambda}_+^{-1}((2r_++d) + \epsilon))  -  {\lambda}_-^{-1}((2r_--d) - \epsilon))\right).
\end{align*}
Denoting 
$$
c:= \left(- {\lambda}_+^{-1}((2r_++d) + \epsilon))  +  {\lambda}_-^{-1}((2r_--d) - \epsilon))\right) >0,
$$
we get
\begin{align*}
\mathsf E_{x,z} (Y_{(k+1)\wedge  N}^2 | {\mathcal F}^Y_{{k}}) - \mathsf E_{x,z} (Y_{k\wedge  N}^2  | {\mathcal F}^Y_{k})
  %\nonumber \\ \nonumber \\ 
\le  -  1(N>k)  \times c.
\end{align*}
Because of the markovian property of $(Y_n)$,
this is equivalent to the second inequality in (\ref{lem11-1}),  as required. 
Summing up from $k=0$ to $k=K$, and taking expectations, we obtain
\begin{align}\label{Y0z0}
c \E\sum_{k=0}^K  1(N>k) \le \E_{x,z}Y_0^2 - \E Y_{K+1}^2,
\end{align}
or, in other words, 
\begin{align*}
c \E (N \wedge (K+1)) \le x^2 - \E Y_{K+1}^2 \le x^2.
\end{align*}
Hence, in the case of $|x| > M_1$ the desired bound (\ref{eN}) follows from the monotone convergence theorem as $K\uparrow \infty$. If $|x|\le M_1$, then $N=0=\E_{x,z}N$.
%Lemma \ref{lem11}  is proved. 

\medskip

Consider now the case $z=1$. Then $T_0>0$ a.s. and $Y_0 = X_{T_0}$. Still, the bound (\ref{Y0z0}) is valid, but now it is unlikely that $Y^2_0=x^2$. We estimate, 
\begin{align*}
& \E_{x,1}Y_0^2 = \E_{x,1}X_{T_0}^2 = x^2 + \E_{x,z}\int_0^{T_0}(2X_tb(X_t) + d)dt  
 \\\\
&\le   x^2 + C \E_{x,z}\int_0^{T_0}1\,dt = x^2 + C\lambda_+^{-1}. 
\end{align*}
So, by virtue of (\ref{Y0z0}) we get 
$$
c\E_{x,1} N \le x^2+C. 
$$
Hence, we obtain (\ref{eN}) in both cases, $z=0,1$. Lemma \ref{lem11} is proved.
\hfill QED

~

\noindent
{\em Proof of theorem \ref{thm2}.} 
It is necessary to estimate the value $\E_{x,z}N^2$ for each $x,z$.
Let us multiply both sides of  the inequality (\ref{lem11-1}) by $n+2$ and use once the inequality 
$n+2 >n+1$:
\begin{equation}\label{Y2ineq2}
1(n+1<N) (n+2)
\E_{Y_{n}} |Y_{(n+1)}|^2 \le 1(n<N) (n+2) |Y_{n}|^2 - c  1(n<N) (n+1),
\end{equation}
or, equivalently, 
$$
 c  1(n<N) (n+1) \le 
1(n<N) (n+2) |Y_{n}|^2 - 1(n<N) (n+2)\E_{Y_{n}} |Y_{n+1}|^2.
$$
So, we estimate, 
\begin{align*}
1(n<N) (n+2) |Y_{n}|^2 = 1(n<N) (n+1) |Y_{n}|^2 (1+\frac1{n+1}),
\end{align*}
so that for $\epsilon <c$ and denoting $\kappa_n := \E1(n<N) (n+1) |Y_{n}|^2$, we obtain
\begin{align*}
& c  E1(n<N) (n+1) \le \kappa_n - \kappa_{n+1} + 
E1(n<N) (n+1) |Y_{n}|^2 \frac1{n+1}
 \\\\
&\le   \kappa_n - \kappa_{n+1} + 
E1(n<N) (n+1) |Y_{n}|^2 \frac1{n+1} 1(|Y_n|^2 \le \epsilon (n+1))
 \\\\
&+ 
E1(n<N) (n+1) |Y_{n}|^2 \frac1{n+1} 1(|Y_n|^2 > \epsilon (n+1)) 
 \\\\
&\le   \kappa_n - \kappa_{n+1} + 
\epsilon E1(n<N) (n+1) \underbrace{1(|Y_n|^2 \le \epsilon (n+1))}_{\le 1}
 \\\\
&+ 
E1(n<N) (n+1) |Y_{n}|^2 \frac1{n+1} \frac{|Y_n|^4}{\epsilon^2 (n+1)^2};
\end{align*}
thus, 
\begin{align*}
 (c-\epsilon)  E1(n<N) (n+1) \le \kappa_n - \kappa_{n+1} 
 %\\\\
+ 
E1(n<N)   \frac{|Y_n|^6}{\epsilon^2 (n+1)^2};
\end{align*}
therefore, with chosen $c$ and $\epsilon$ we obtain
\begin{align*}
& (c-\epsilon)  E \sum 1(n<N) (n+1) \le \kappa_0 - \kappa_{N} 
 %\\\\
+ 
E\sum_n1(n<N)   \frac{|Y_n|^6}{\epsilon^2 (n+1)^2}
 \\\\
&\le  \kappa_0 - \kappa_{N} 
 %\\\\
+ 
E\sum_n   \frac{C(|x|^6+1)}{\epsilon^2 (n+1)^2}
\le x^2 + C(x^6+1).
\end{align*}
In other words, 
\begin{align*}
 (c-\epsilon)  E N^2/2 \le 
 (c-\epsilon)  E N(N+1)/2 
\le x^2 + C(x^6+1).
\end{align*}

\noindent
By virtue of (\ref{etau2}) and (\ref{eN}), and since $x^2\le x^6+1$ for any $x$,  this imples the desired bound (\ref{e3a}). Theorem \ref{thm2} is proved 
\hfill QED

\ifhide
\medskip

{\color{blue}
Consider the case $z=1$. Since $N\ge 0$, we estimate using the previous and the strong Markov property, 
\begin{align*}
\E_{x,1} N^2 
\le C\E_{x,1}(1+X_{T_0}^6).
\end{align*}
Further, 
\begin{align*}
&\E_{x,1} X_{T_0}^6 - x^6 = \E_{x,1}\int_0^{T_0} (6X_t^5b(X_t) + 15 X^4_t)dt 
 \\\\
&=\int_0^{\infty} \E_{x,1}1(t<T_0) (6(X_t^1)^5b(X_t) + 15 (X_t^1)^4)dt 
 \\\\
&=\int_0^{\infty} \E_{x,1}1(t<T_0) \E_{x,1}(6(X_t^1)^5b(X_t) + 15 (X_t^1)^4)dt 
 \\\\
&\le C \int_0^{\infty} \exp(-\lambda_+ t) \E_{x,1}(X_t^1)^4dt
 \\\\
&\stackrel{\text{lemma \ref{lem50}}}
\le 
C \int_0^{\infty} \exp(-\lambda_+ t) (x^4 + Cx^2t + Ct^2 +Ct)dt 
 \\\\
&\le C(x^4+x^2+1) \le C(x^6+1). 
\end{align*}
Hence, theorem follows in the case $z=1$ as well. \hfill QED 
}
\fi

\ifhide
\begin{Remark}
For the $k$-moment of $\tau$ and, respectively, $2k$-moment for $Y_n$ we would need to use  the condition
$$%\begin{equation}\label{c2}
\frac{(2k r_- - k(2k-1)d)}{\lambda_-} > \frac{(2kr_++k(2k+1)d)}{\lambda_1}, 
$$%\end{equation}
or, equivalently, 
$$%\begin{equation}\label{c2}
\frac{(2 r_- - (2k-1)d)}{\lambda_-} > \frac{(2r_++(2k+1)d)}{\lambda_1}, 
$$%\end{equation}
also combined with (\ref{c2R}) and (\ref{c2R2}).
\end{Remark}
\fi

%\section*{Acknowledgments}
%\noindent
%This research was funded by Russian Foundation for Basic Research grant 20-01-00575a.

%\medskip

%\noindent
%{\bf Data availability statement. } Data sharing is not applicable to this article as no datasets were generated or analyzed during the current study.

%
% ---- Bibliography ----
%

\end{document}